\documentclass[english,a4,twoside]{scrartcl}

\usepackage{float}
\usepackage{graphicx}
\usepackage{tikz}
\usepackage{pgflibraryarrows}

\usepackage[margin=20pt]{caption}

\usepackage{amsmath}
\usepackage{amssymb}
\usepackage{amsthm}

\usepackage[english]{babel}
\usepackage{flafter}
\usepackage{array}
\usepackage{paralist}

\graphicspath{{graphics/}}

\newtheoremstyle{break}%
{9pt}{9pt}%
{\itshape\raggedright}%
{0pt}%
{\sffamily\bfseries}{}%
{\newline}%
{}%

\theoremstyle{break}

\newtheorem{theorem}{Theorem}

\newtheorem*{corollary*}{Corollary}
\newtheorem{definition}[theorem]{Definition}
\newtheorem*{definition*}{Definition}

\newtheorem{proposition}[theorem]{Proposition}
\newtheorem*{proposition*}{Proposition}
\newtheorem{property}[theorem]{Property}
\newtheorem*{property*}{Property}
\newtheorem{remark}[theorem]{Remark}
\newtheorem*{remark*}{Remark}

\newcommand\qedheretext[1]{\qedhere\raisebox{-0.5em}{\;#1}}

\newcommand\MST{\operatorname{MST}}
\newcommand{\set}[1]{\left\{ {#1} \right\}}
\newcommand{\xbinom}{\tbinom}

\newcommand{\ignore}[1]{\relax}
\newcommand{\order}[1]{\langle #1 \rangle}

\title{A tight bound on the collection of edges in MSTs of induced subgraphs}

\author{Gregory B. Sorkin, Angelika Steger and Rico Zenklusen}

\begin{document}

\maketitle

\begin{abstract}
Let $G=(V,E)$ be a complete $n$-vertex graph with distinct positive edge
weights.  We prove that for $k\in\{1,2,\dots,n-1\}$, the set consisting of the
edges of all minimum spanning trees (MSTs) over induced subgraphs of $G$ with
$n-k+1$ vertices has at most $nk-\tbinom{k+1}{2}$ elements. This proves a
conjecture of Goemans and Vondrak~\cite{goemans_2005_covering}. We also show
that the result is a generalization of Mader's Theorem, which bounds the
number of edges in any edge-minimal $k$-connected graph.
\end{abstract}

\section{Introduction}

Let $G=(V,E)$ be a complete $n$-vertex graph with distinct
positive edge weights.  
For any set $X\subseteq V$,
denote by $G[V \setminus X]$ the subgraph of $G$ induced by~$V \setminus X$.
We will also sometimes write this graph as $(V \setminus X, E)$, 
ignoring edges in $E$ incident on vertices in~$X$.
$\MST(G[V \setminus X])$ denotes the set of edges in the graph's 
minimum spanning tree.
(The MST is unique due to the assumption that the edge
weights are distinct.)

For $k\in\{1,2,\dots,n-1\}$, define
\begin{equation*}
M_k(G)=\bigcup_{X\subseteq V, \: |X|=k-1}\MST(G[V\setminus X])\;.
\end{equation*}

Note that for $k=1$ we have $M_1(G)=\MST(G)$. In
\cite{goemans_2005_covering}, Goemans and Vondrak considered 
the problem of finding a
sparse set of edges which, with high probability,
contain the MST of a random subgraph of~$G$. In this context
they proved an upper bound on $M_k(G)$, namely that 
$|M_k(G)|<(1+\frac{e}{2})kn$, and they conjectured that one
should be able to improve the bound to $|M_k(G)|\leq nk
-\tbinom{k+1}{2}$. In this paper we prove this conjecture.

\begin{theorem}\label{thm:original_conjecture}
For any complete graph $G$ on $n$ vertices with distinct positive edge weights,
\begin{equation}
|M_k(G)|\leq nk-\xbinom{k+1}{2} . \label{maineq}
\end{equation}
\end{theorem}

As Goemans and Vondrak recognized, the bound is tight:
for any $n$ and $k$ it is easy to produce edge weights 
giving equality in~\eqref{maineq}.
One way is to fix an arbitrary
set $V'\subseteq V$ with cardinality~$k$, and partition the
edges $E$ into three sets $E_0$, $E_1$ and $E_2$ where,
for $i\in\{0,1,2\}$, $E_i$ contains all edges of $E$ having
exactly $i$ endpoints in~$V'$. Assign arbitrary distinct positive
weights to the edges in $E$
such that all weights on $E_2$ are smaller than those
on~$E_1$, which in turn are smaller than those on~$E_0$. It
can easily be verified that $M_k(G)=E_2\cup E_1$ and thus
$|M_k(G)|=nk-\tbinom{k+1}{2}$.

Theorem~\ref{thm:original_conjecture}'s assumption
that $G$ is complete is not meaningfully restrictive.
If $G$ is such that deletion of some $k-1$ vertices leaves it 
disconnected, then the notion of $M_k(G)$ does not make sense;
otherwise, it does not matter if other edges of $G$ are 
simply very costly or are absent.

The bound of Theorem~\ref{thm:original_conjecture} applies equally 
if we consider the edge set of MSTs of induced subgraphs of size 
\emph{at most} $n-k+1$ (rather than exactly that number).
This is an immediate consequence of the following remark.

\begin{remark}
For any complete graph $G$ on $n$ vertices with distinct positive edge weights,
and $k \in \set{1,2,\dots,n-2}$,
$M_{k+1}(G) \supseteq M_{k}(G)$.
\end{remark}
\begin{proof}
We will show that any edge $e$ in $M_{k}(G)$ is also in $M_{k+1}(G)$.
By definition, $e \in M_{k}(G)$ means that there is 
some vertex set~$X$ of cardinality $|X|=k-1$ for which 
$e \in \MST(G_k)$, where $G_k = G[V \setminus X]$.

Consider any leaf vertex $v$ of~$\MST(G_k)$, with neighbor~$u$.
We claim that deleting $v$ from $G_k$ 
(call the resulting graph $G_{k+1}$) results in 
the same MST less the edge $\set{u,v}$, i.e.,
that $\MST(G_{k+1}) = \MST(G_k) \setminus \set{\set{u,v}}$.
This follows from considering the progress of Kruskal's algorithm 
on the two graphs.
Before edge $\set{u,v}$ is added to $\MST(G_k)$, 
the two processes progress identically:
every edge added to $\MST(G_k)$ is also a cheapest edge 
for the smaller graph~$G_{k+1}$.
The edge~$e$, added to $\MST(G_k)$, of course has no parallel in~$G_{k+1}$. 
As further edges are considered in order of increasing cost,
again, every edge added to $\MST(G_k)$ will also be added to $\MST(G_{k+1})$, 
using the fact that none of these edges is incident on~$v$.

Thus, if $v$ is not a vertex of $e$, then $e \in \MST(G_{k+1})$. 
Since $\MST(G_k)$ has at least two leaves, 
it has at least one leaf $v$ not in~$e$, unless $\MST(G_k)=e$, 
which is impossible since $G_k$ has at least 3 vertices.
\end{proof}

\subsection*{Outline of the paper}

In Section~\ref{sec:k_constructible} we define a ``$k$-constructible'' graph,
and show that
every graph $(V,M_k(G))$ is $k$-constructible,
and every $k$-constructible graph is a subgraph of some graph $(V,M_k(G))$.
This allows a simpler reformulation of Theorem~\ref{thm:original_conjecture}
as Theorem~\ref{thm:edge_bounding},
which also generalizes a theorem of Mader \cite{mader_1971_minimale}.
We prove Theorem~\ref{thm:edge_bounding} in Section~\ref{sec:mainproof}.

\section{$k$-constructible graphs}\label{sec:k_constructible}

We begin by recalling Menger's theorem for undirected
graphs, which motivates our definition of $k$-constructible
graphs. Two vertices in an undirected graph are called
$k$-connected if there are $k$ (internally) vertex-disjoint paths
connecting them.

\begin{theorem}[Menger's theorem]
Let $s,t$ be two vertices in an undirected graph $G=(V,E)$ such that
$\{s,t\}\not\in E$. Then $s$
and $t$ are $k$-connected in $G$ if and only if after
deleting any $k-1$ vertices (distinct from $s$ and $t$), $s$ and
$t$ are still connected.
\end{theorem}

\begin{definition}[$k$-constructible graph]
A graph $G=(V,E)$ is called $k$-constructible if there exists
an ordering $O=\order{e_1,e_2,\dots,e_m}$ of the edges in $E$
such that for all $i\in \{1,2,\dots,m\}$ the graph
$(V,\{e_1,e_2,\dots,e_{i-1}\})$ contains at most $k-1$ vertex-disjoint paths 
between the two endpoints of~$e_i$. We say
that $O$ is a $k$-construction order for the graph $G$.
\end{definition}

Note that $1$-constructible graphs are forests, and
edge-maximal $1$-constructible graphs are spanning trees. We therefore
have in particular that graphs of the form $M_1(G)$ 
(i.e., MSTs, recalling the $G$ is complete) 
are edge-maximal $1$-constructible graphs. 
A slightly weaker statement 
is true for all~$k$: 
every graph $M_k(G)$ is $k$-constructible 
(Theorem~\ref{thm:equivalence}.\ref{5i}), 
and every $k$-constructible graph is a subgraph of some graph~$M_k(G)$
(Theorem~\ref{thm:equivalence}.\ref{5ii}).

Note that a stronger statement, that the graphs of the form $M_k(G)$ 
are exactly the edge-maximal $k$-constructible graphs, is not true. 
To see this consider a cycle $C_4$ of length four. Assign weights
$1,\dots,4$ to these four edges (in arbitrary order) and weights $5,6$
to the remaining edges of the complete graph on  four vertices. It
is easily checked that $M_2(G)=C_4$. But $M_2(G)$ is not edge-maximal,
as a diagonal to the cycle $C_4$ can be added without destroying
$2$-constructibility.
   
\begin{theorem} \label{thm:equivalence}
\hspace*{0mm}\vspace*{-\baselineskip}
\begin{compactenum}[i)]
\item For every complete graph $G=(V,E)$ with distinct
  positive edge weights, $(V,M_k(G))$
  is $k$-constructible.
  \label{5i}
\item Let $G=(V,E)$ be $k$-constructible. Then there exist distinct
  positive edge weights 
  for the complete graph 
  $\widetilde{G} = (V,\widetilde{E})$ such that
  $E\subseteq M_k(\widetilde{G})$.
  \label{5ii}
\end{compactenum}
\end{theorem}

\begin{proof}
Part~\eqref{5i}:
Let $G=(V,E)$ be a complete graph on $n$ vertices with
distinct positive edge weights. 
Let $\order{e_1,e_2,\dots,e_{\tbinom{n}{2}}}$ be the
ordering of the edges in $E$ by increasing edge weights
and $O=\order{e_{r_1},e_{r_2},\dots,e_{r_{|M_k(G)|}}}$ be
the ordering of the edges in $M_k(G)$ by increasing edge
weights. We will now show that $O$ is a
$k$-construction order for $(V,M_k(G))$. Let
$i\in\{1,2,\dots,|M_k(G)|\}$. As $e_{r_i} \in M_k(G)$ there exists a
set $X\subseteq V$ with $|X|=k-1$ and
$e_{r_i} \in \MST(G\setminus X)$, implying that the two endpoints
of $e_{r_i}$ are not connected in the graph
$(V\setminus X, \{e_1,e_2,\dots,e_{r_i-1}\})$. 
By Menger's
theorem, this implies that there are at most $k-1$ vertex-disjoint paths 
between the two endpoints of $e_{r_i}$
in $(V,\{e_1,e_2,\dots,e_{r_i-1}\})$. This statement remains
thus true for the subgraph
$(V,\{e_{r_1},e_{r_2},\dots,e_{r_{i-1}}\})$.
The ordering $O$ is thus a $k$-construction order for $(V,M_k(G))$.

Part~\eqref{5ii}:
Conversely let $G=(V,E)$ be a $k$-constructible graph with
$k$-construction order $O=\order{e_1,e_2,\dots,e_{|E|}}$.
Let $(V,\widetilde{E})$ be the complete graph on~$V$. 
We assign the following edge weights $\widetilde{w}$ to the edges
in~$\widetilde{E}$. We assign the weight $1$ to $e_1$, $2$ to $e_2$
and so on. The remaining edges $\widetilde{E}\setminus E$
get arbitrary distinct weights greater than~$|E|$. In order to show
that the graph $\widetilde{G}=(V,\widetilde{E},\widetilde{w})$
satisfies $E\subseteq M_k(\widetilde{G})$ consider an arbitrary edge
$e_i\in E$ and let $C\subseteq V$ with
$|C|=k-1$ be a vertex set separating the two endpoints of
$e_i$ in the graph $G_{i-1}=(V,\{e_1,e_2,\dots,e_{i-1}\})$.
Applying Kruskal's algorithm to
$\widetilde{G}[V\setminus C]$, 
the set of all edges considered before $e_i$ is contained 
in $E(G_{i-1})$, leaving the endpoints of $e_i$ separated, 
so $e_i$ will be accepted: 
$e_i \in \MST(\widetilde{G}[V\setminus C]) \subseteq M_k(\widetilde G)$.
\end{proof}

We remark that the first part of the 
foregoing proof shows an efficient construction of~$M_k(G)$: 
follow a generalization of Kruskal's algorithm, 
considering edges in order of increasing weight, 
adding an edge if (prior to addition) its endpoints are 
at most $(k-1)$-connected.
Connectivity can be tested as a flow condition, so that the algorithm runs 
in polynomial time --- far more efficient than the naive 
$\Omega\left( \tbinom n k \right)$ protocol suggested by the definition
of~$M_k(G)$.
This again was already observed in~\cite{goemans_2005_covering}.

By Theorem~\ref{thm:equivalence},
the following theorem
is equivalent to Theorem~\ref{thm:original_conjecture}.

\begin{theorem}\label{thm:edge_bounding}
For $k\geq 1$, every $k$-constructible graph $G=(V,E)$
with $n \geq k+1$ vertices satisfies
\begin{equation}
|E|\leq nk-\xbinom{k+1}{2}\;.  \label{equivalent}
\end{equation}
\end{theorem}

Theorem~\ref{thm:edge_bounding} generalizes a result of 
Mader~\cite{mader_1971_minimale}, 
based on results in~\cite{halin_1969_theorem}, concerning 
``$k$-minimal'' graphs 
(edge-minimal $k$-connected graphs).
Every $k$-minimal graph is $k$-constructible, since
every order of its edges is a $k$-construction order.
The following theorem is thus a corollary of Theorem~\ref{thm:edge_bounding}.

\begin{theorem}[Mader's theorem] \label{Mader}
Every $k$-minimal graph with $n$ vertices 
has at most $nk-\tbinom{k+1}{2}$ edges.
\end{theorem}

Note that Mader's theorem (Theorem~\ref{Mader}) is 
weaker than Theorem~\ref{thm:edge_bounding}, because
while every $k$-minimal graph is $k$-constructible, the converse is false:
not every $k$-constructible graph is $k$-minimal. 
An example with $k=2$ is a cycle $C_4$ with length four with an
additional diagonal~$e$. The vertex set remains $2$-connected even upon
deletion of the edge $e$, so the graph is not $2$-minimal, but it
is $2$-constructible (by any order where $e$ is not last).

\section{Proof of the main theorem} \label{sec:mainproof}

In this section we prove Theorem~\ref{thm:edge_bounding}.
We fix $k$ and prove the theorem by induction on~$n$. 
The theorem is trivially true for
$n = k+1$, so assume that $n \geq k+2$
and that the theorem is true for all smaller values of~$n$.
We prove \eqref{equivalent} for a
$k$-constructible graph $G=(V,E)$ on $n$
vertices and $m$ edges which, without loss of generality, 
we may assume is edge-maximal 
(no edges may be added to $G$ leaving it $k$-constructible).
Fix a $k$-construction order
\begin{align*}
O &= \order{e_1,e_2,\dots,e_m}
\end{align*}
of $G$ and (for any $i \leq m$) let
$G_i=(V,\{e_1,e_2,\dots,e_i\})$. 
Also fix a set
$C\subseteq V$ of size $|C|=k-1$ such that the two
endpoints of $e_m$ lie in two different components $Q^1,Q^2
\subseteq V$ of $G_{m-1}[V\setminus C]$
(the set $C$ exists by $k$-constructibility of $G$ and
Menger's theorem). The edge maximality of $G$
implies that $Q^1,Q^2,C$ form a
partition of~$V$. Let $V^1=Q^1\cup C$ and $V^2=Q^2\cup C$.
(If there were a third component $Q^3$ then, even after adding~$e_m$,
any $v_1 \in Q^1$ and $v_3 \in Q^3$ are at most $(k-1)$-connected
and so the edge $\set{v_1,v_3}$ could be added,
contradicting maximality.)

Our goal is to define two graphs 
$G^1=(V^1,E^1)$ and $G^2=(V^2,E^2)$ that satisfy
the following property.
\begin{property}%
\label{prope:requested_E1_E2}\hspace*{0mm}%
\vspace*{-\baselineskip}
\begin{compactitem}
\item $G^1$ and $G^2$ are both $k$-constructible.
\item $E^1$ contains all edges of $G[V^1]$.
\item $E^2$ contains all edges of $G[V^2]$.
\item For every pair of vertices $c_1,c_2
\in C$ not connected by an edge in~$G$, there is an
edge $\{c_1,c_2\}$ in either $E^1$ or in $E^2$ (but not both).
\end{compactitem}
\end{property}

If we can find graphs $G^1$ and $G^2$
satisfying Property~\ref{prope:requested_E1_E2}, then
the proof can be finished as follows. Note that we have the
following equality:
\begin{equation*}
|E^1|+|E^2| = (m-1)+|G[C]| + \left(\xbinom{k-1}{2}-|G[C]|\right)\;.
\end{equation*}
The term $m-1$ comes from the fact that $E^1\cup E^2$ covers
all edges of $G$ except~$e_m$, the term $|G[C]|$ represents the
double counting of edges contained in~$C$, and the last term
counts the edges 
which are covered by $E^1$ and $E^2$ but not in~$G$.

We therefore have
\begin{equation*}
m=1+|E^1|+|E^2|-\xbinom{k-1}{2}\;.
\end{equation*}

Applying the inductive hypothesis on $G^1$ and $G^2$ 
(which by Property~\ref{prope:requested_E1_E2} are $k$-constructible)
we get the desired result:

\begin{align*}
m&\leq 1+\left(|V^1|k-\xbinom{k+1}{2}\right)
+\left(|V^2|k-\xbinom{k+1}{2}\right)-\xbinom{k-1}{2}\\
&\leq 1+(n+k-1)k-2\xbinom{k+1}{2}-\xbinom{k-1}{2}
\\
&=nk-\xbinom{k+1}{2}\;,
\end{align*}
where in the second inequality we have used $|V_1|+|V_2|=n+|C| = n+k-1$.

We will finally concentrate on finding 
$G^1=(V^1,E^1)$ and
$G^2=(V^2,E^2)$ satisfying Property~\ref{prope:requested_E1_E2}.

Let $B=\tbinom{C}{2} \setminus E$ be the set of all 
anti-edges in~$G[C]$.
($\tbinom{C}{2}$ denotes the set of unordered pairs of elements of~$C$.)
For $\set{c_1,c_2} \in B$, let $\ell(c_1,c_2)$ be the smallest value of $i$ 
such that $c_1$ and $c_2$ are $k$-connected in~$G_i$. 
(Considering $k$ vertex-disjoint paths between $c_1$ and $c_2$ in~$G_i$, 
and noting that deletion of the single edge $e_i$ leaves them at least
$k-1$ connected, it follows that $c_1$ and $c_2$ 
are precisely $(k-1)$-connected in~$G_{i-1}$.)
Define $B_i = \set{ \set{c_1,c_2} \colon \ell(c_1,c_2)=i }$. 
Since by edge maximality of $G$ \emph{every} pair $\set{c_1,c_2}$ 
is $k$-connected in $G_m=G$, 
it follows that $B_1,B_2,\dots,B_m$ form a partition of~$B$.

Our basic strategy to define the graphs $G^1$ and $G^2$ 
(and appropriate orderings of their edges
which prove that they are $k$-constructible) is as follows.
In a particular way, 
we will partition each $B_i$ as $B_i^1 \cup B_i^2$, 
and determine orders $O_i^1$ and $O_i^2$ on their respective edges. 
Let $G^1$ be the graph constructed by the order
\begin{align} \label{O1}
O^1 &= \order{ e_1, O_1^1, e_2, O_2^1, \ldots, e_m, O_m^1 } ,
\end{align}
where (recalling that $G^1$ has vertex set $V^1$) 
we ignore any edge $e_i \notin \tbinom{V^1}2$.
(There is no issue with edges from $O_i^1$, as these
belong to $\binom C2 \subseteq \binom{V^1}2$.) 
Define $G^2$ symmetrically.
We need to show that the graphs $G^1$ and $G^2$ satisfy
Property~\ref{prope:requested_E1_E2};
the central point will be to ensure that 
$O^1$ is a $k$-construction order for $G^1$,
and $O^2$ for $G^2$.
(By definition of the edges $B_i$, note that every edge $e \in O_i^1$ 
when added after $e_i$ in the order $O$ violates $k$-constructibility, 
but in the following we show how $O_i^1,O_i^2$ can be
chosen such that it will not violate $k$-constructibility in $G^1$; 
likewise for edges $e \in O_i^2$ and $G^2$.)

To show that $O^1$ and $O^2$ are $k$-construction orders 
we need to check that, just before an edge is added,
its endpoints are at most $(k-1)$-connected.
To prove this, we distinguish between edges $e_i \in E$ and edges $e \in B$. 
We first dispense with the easier case of an edge $e_i \in E$.
Proposition~\ref{propo:B_invariance} shows that
(for \emph{any} orders $O_i$ of $B_i$)
in the edge sequence $\order{ e_1,O_1,\ldots,e_m,O_m }$,
every edge $e_i$ has endpoints which are at most $(k-1)$-connected
upon its addition to the graph
$(V, \set{e_1,O_1,\ldots,e_{i-1},O_{i-1}})$. 
It follows that the endpoints are also at most $(k-1)$-connected 
upon the edge's addition to $G^1$ (respectively,~$G^2$), i.e., in the graph
$(V^1, \set{e_1,O_1^1,\ldots,e_{i-1},O_{i-1}^1})$,
where as usual we disregard edges not in $\binom{V_1}2$.

\begin{proposition}\label{propo:B_invariance}
Let $i\in\{1,2,\dots,m\}$ and $v_1,v_2 \in V$
such that $\{v_1,v_2\}$ is not an edge in $G_{i-1}$.
If the maximum number of vertex-disjoint paths between
$v_1$ and $v_2$ in $G_{i-1}$ is $r\leq k-1$, then the maximum
number of vertex-disjoint paths between $v_1$ and $v_2$ in
the graph 
$(V,\{e_1,e_2,\dots,e_{i-1}\} \cup \bigcup_{l=1}^{i-1} B_l)$ is~$r$, too.
\end{proposition}

\begin{proof}
For any $i,v_1,v_2$ as above,
let $S\subseteq V$, $|S|=r$, be a set separating $v_1$ and
$v_2$ in~$G_{i-1}$. As $|S|=r<k$, $S$ cannot separate
two $k$-connected vertices in~$G_i$. This implies that any two
vertices in $V\setminus S$ that are $k$-connected in $G_{i-1}$
lie in the same connected component of \mbox{$G_{i-1}[V\setminus S]$}. As
every edge in $\bigcup_{l=1}^{i-1} B_l$ connects two vertices
that are $k$-connected in~$G_{i-1}$, adding the
edges $\bigcup_{l=1}^{i-1} B_l$ to \mbox{$G_{i-1}[V\setminus S]$} 
\emph{does not change the component structure} of $G_{i-1}[V\setminus S]$. The
set $S$ thus remains a separating set for $v_1$ and $v_2$ in the graph
$(V,\{e_1,e_2,\dots,e_{i-1}\}\cup\bigcup_{l=1}^{i-1} B_l)$,
proving that $v_1$ and $v_2$ are at most $r$-connected in this graph.\\
\qedheretext{Proposition \ref{propo:B_invariance}}
\end{proof}

With Proposition~\ref{propo:B_invariance} addressing edges $e_i \in E$, 
to ensure $k$-constructibility of $O^1$ and $O^2$,
it suffices 
to choose for $j\in\{1,2\}$ and $i\in\{1,2,\dots,m\}$ the
orders $O_i^j$ in such a way that successively
adding any edge $e \in O_i^j$ to the graph $G_i[V^j]$ 
connects two vertices which were at most $(k-1)$-connected.

\begin{figure}[H]
\begin{center}
\includegraphics{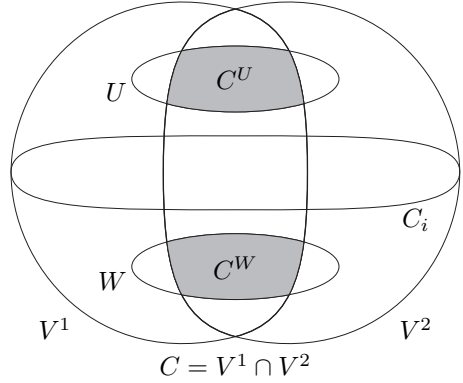}
\end{center}\vspace*{-\baselineskip}
\caption{Sets defined to prove 
Propositions \ref{propo:bipartite_Bi}--\ref{propo:next_vertex}.%
}
\label{fig:CU_CW}
\end{figure}

Let $C_i \subseteq V$ with $|C_i|=k-1$ a set separating the 
endpoints of $e_i$ in the graph $G_{i-1}$. 
Let $U,W \subseteq V$ be the two components
of $G_{i-1}[V\setminus C_i]$ containing
the two endpoints of the edge~$e_i$. 
We define $C^U=C\cap U$, $C^W=C\cap W$. 
Figure~\ref{fig:CU_CW} illustrates these sets.

The following proposition shows that the edges $B_i$ form a bipartite graph.

\begin{proposition}\label{propo:bipartite_Bi}
\begin{equation*}
B_i \subseteq C^U \times C^W
\end{equation*}
\end{proposition}

\begin{proof}
Suppose by way of contradiction that 
$\exists e \in B_i \setminus (C^U \times C^W)$.
Let 
\begin{align*}
O' &= \order{ e_1, \ldots, e_{i-1}, e, e_i, \ldots, e_m },
\end{align*}
the edge order obtained by inserting $e$
immediately before $e_i$ in the original order 
$O=\order{ e_1,e_2,\dots,e_m }$. 
We will show that $O'$ is a $k$-construction order, thus contradicting
the edge maximality of~$G$. 
For edges up to $e_{i-1}$ this is immediate from the fact that $O$
is a $k$-construction order. Proposition~\ref{propo:B_invariance}
shows that edges $e_{i+1}$ and later do not violate $k$-constructibility.
(Literally, Proposition~\ref{propo:B_invariance} applies to the order 
$\order{e_1,\ldots, e_i, e, e_{i+1}, \ldots, e_m}$
rather than to $O'$, but for edges $e_{i+1}$ and later the swap of 
$e_i$ and $e$ is irrelevant.)
The edge $e$ itself does not violate $k$-constructibility,
since by the definition of $B_i$ its two endpoints are at most $k-1$
connected in~$G_{i-1}$. 
This leaves only edge $e_i$ to check, 
but since $e \notin U \times W$, 
$C_i$ remains a separating set with cardinality $k-1$
for the two endpoints of $e_i$ in the graph
$(V,\set{e_1,e_2,\dots,e_{i-1},e})$. 
Thus $O'$ is a $k$-construction order, giving the desired contradiction.
\end{proof}

We will now describe a method for constructing the orders $O_i^1$, $O_i^2$.
Our approach is to define an order $L=\order{ v_1,v_2,\dots,v_r }$ 
on (a subset of) the \emph{vertices} of $C^U\cup C^W$ 
and to assign to every vertex
$v\in C^U\cup C^W$ a label $\alpha(v)\in\{1,2\}$.
The two orders $O_i^1$, $O_i^2$ are then defined as follows. 
We begin with $O_i^1,O_i^2=\emptyset$ and add all edges in $B_i$
which are incident to $v_1$ at the end of $O_i^{\alpha(v_1)}$ in any order. 
In the next step all edges of $B_i$ which are incident to $v_2$ and
not already assigned to one of the orders $O_i^1,O_i^2$ are
added at the end of $O_i^{\alpha(v_2)}$ in any order.
This is repeated until all edges are assigned.

In what follows we show how to 
choose a vertex order $L$ and labels $\alpha$ so that 
$O^1$ and $O^2$ are $k$-construction orders.
Just as $O^1$ and $O^2$ are built iteratively, so is~$L$, 
starting with $L=\emptyset$.

For any $X\subseteq C^U\cup C^W$, we define $B_i(X)$ to be the set of edges 
in $B_i$ incident on vertices in $X$, i.e.,
$B_i(X) = \{e\in B_i\mid e\cap X \neq \emptyset\} $.

\begin{proposition}\label{propo:edge_adding_freedom}
Let $j\in\{1,2\}$ and $X \subseteq C^U\cup C^W$. We then
have that $\forall e\in B_i\setminus B_i(X)$ there
are at most $|C_i\cap V^j|+|X|$ vertex-disjoint paths
between the two endpoints of $e$ in the graph
$(V,\{e_1,e_2,\dots,e_i\}\cup B_i(X))[V^j]$.
\end{proposition}

\begin{proof}
Observe that the set $(C_i\cap V^j) \cup X$ separates the two
endpoints of the edge $e$ in the graph
$(V,\{e_1,e_2,\dots,e_i\}\cup B_i(X))[V^j]$. As this set has cardinality
$|C_i\cap V^j|+|X|$ the result follows by Menger's theorem.
\end{proof}

Let $X^1$ be the set of vertices labeled 1 contained in the 
partially constructed~$L$, and $X^2$ those labeled~2.
If we can find a vertex $v\in (C^U\cup C^W)\setminus (X^1\cup X^2)$ 
where the number of ``new'' edges incident on $v$ satisfies
\begin{equation}\label{eq:next_vertex}
|B_i(v)\setminus (B_i(X^1 \cup X^2))|\leq
k-1-\min\{|C_i\cap V^1|+|X^1|,|C_i\cap V^2|+|X^2|\}
\end{equation}
then by Proposition~\ref{propo:edge_adding_freedom},
adding $v$ at the end of the current order
$L$ and labeling it 
$\arg\min_{j\in\{1,2\}}\{|C_i\cap V^j|+|X^j|\}$
does not violate $k$-constructibility of the orders $O^1$
and~$O^2$.

The following proposition shows that, until the process is complete
(until $B_i(X^1 \cup X^2) = B_i$), such a vertex $v$ can always be found.

\begin{proposition}\label{propo:next_vertex}
Let $X^1,X^2\subset C^U \cup C^W$ be two disjoint sets. 
If $B_i(X^1 \cup X^2) \subsetneq B_i$, 
then there exists a vertex $v\in (C^U\cup C^W)\setminus (X^1\cup X^2)$
that satisfies~(\ref{eq:next_vertex}).
\end{proposition}
\begin{proof}
Note that $C^U$, $C^W$, and $C \cap C_i$ are disjoint and contained in~$C$, so
\begin{align} 
|C^U|+|C^W|+|C\cap C_i| & \leq |C| = k-1 \; , \label{first}
\end{align}
where $|C|=k-1$ by definition.  Also,
\begin{align}
|V^1\cap C_i|+|V^2\cap C_i|-|C\cap C_i|= |C_i| = k-1 \; . \label{second}
\end{align}
>From the fact that the right side of \eqref{first} is equal to 
$2(k-1)$ minus that of \eqref{second}, we get
\begin{align}
|C^U|+|C^W| &\leq (k-1-|V^1\cap C_i|)+(k-1-|V^2\cap C_i|) \; . \label{third}
\end{align}
By disjointness of $C^U$ and $C^W$,
\begin{align}
|C^U \setminus & (X^1 \cup X^2)| + |C^W \setminus (X^1 \cup X^2)| \label{4}
\\ & =
|C^U| + |C^W| - |X^1| - |X^2| \notag
\\ & \leq
(k-1-|V^1\cap C_i| - |X^1|)+(k-1-|V^2\cap C_i| - |X^2|) \label{5} \; ,
\end{align}
using \eqref{third} in the last inequality.
Thus, the smaller summand in \eqref{4} is at most the larger summand
in~\eqref{5}, and without loss of generality we suppose that 
\begin{equation}
|C^U\setminus (X^1\cup X^2)|\leq k-1-|V^1\cap C_i|-|X^1|\;. \label{6}
\end{equation}
By the hypothesis $B_i(X^1 \cup X^2) \subsetneq B_i$,
there is an edge $e \in B_i \setminus B_i(X^1 \cup X^2)$; 
by Proposition~\ref{propo:bipartite_Bi}, $e=\{u,w\}$
with $u \in C^U$ and $w \in C^W$;
and by definition of $B_i(X^1 \cup X^2)$,
$u,w \notin X^1 \cup X^2$,
i.e., $u \in C^U \setminus (X^1 \cup X^2)$
and $w \in C^W \setminus (X^1 \cup X^2)$.
Then $v=w$ satisfies \eqref{eq:next_vertex}
because the new edges on $w$ must go to so-far-unused vertices in $C^U$:
\begin{align*}
 | B_i(w) \setminus B_i(X^1 \cup X^2) |
  &\leq | C^U \setminus (X^1 \cup X^2) | \; ,
\end{align*}
whence \eqref{6} closes the argument.
\end{proof}

Therefore there always exist two $k$-construction
orders $O^1,O^2$ as desired, which completes the proof of 
Theorem~\ref{thm:edge_bounding}.

\section*{Acknowledgment}
The authors are grateful to Michel Goemans for bringing the problem 
to their attention.

\bibliographystyle{plain}
\bibliography{mst_edge_bound}

\end{document}